\documentclass[11pt]{article}
\usepackage{amsmath,amscd,amssymb,amsfonts,color,epsfig,amsthm,wrapfig}

\date{21 October 2003}

\title{%
The Satisfiability Threshold of \\ Random 3-SAT Is at Least 3.52
}
\author{%
{\sl  MohammadTaghi Hajiaghayi\/}%
\thanks{MIT Laboratory for Computer Science,
200 Technology Square,
Cambridge MA 02139, U.S.A.
\texttt{hajiagha@mit.edu}}
\medskip \and
{\sl Gregory B. Sorkin\/}%
\thanks{Department of Mathematical Sciences, IBM T.J. Watson Research
Center, Yorktown Heights NY 10598, U.S.A.
\texttt{sorkin@watson.ibm.com}}
}

\newcommand{\xxxx}[1]{}

\let\realuparrow=\uparrow
\let\realdownarrow=\downarrow
\let\realupdownarrow=\updownarrow
\def\arrowspacing{\hskip -0.25em}
\def\uparrow{\realuparrow\arrowspacing}
\def\downarrow{\realdownarrow\arrowspacing}
\def\updownarrow{\realupdownarrow\arrowspacing}

\renewcommand{\H}{\end{document}}

\renewcommand{\leq}{\leqslant}
\renewcommand{\geq}{\geqslant}

\newcommand{\ignore}[1]{}

\def\bol{Bollob\'{a}s}

\begin{document}
\maketitle

\begin{abstract}
We prove that a random 3-SAT instance with clause-to-variable density
less than 3.52 is satisfiable with high probability.
The proof comes through an algorithm which selects (and sets) a variable
depending on its degree and that of its complement.
\end{abstract}

\vspace{-0.25cm}

\section{Introduction}
There is much interest in understanding ``phase transitions''
in mathematics, computer science, and mathematical physics,
and in particular the $k$-SAT phase transition.
In the standard model for random $k$-SAT,
a random $k$-CNF formula $F(n,cn)$
with $n$ variables and density $c$
has $m=cn$ random clauses
independently selected uniformly at random,
with replacement, from among the $2^k \binom{n}{k}$ proper clauses of
length $k$.
The Satisfiability Threshold Conjecture asserts that
for each $k\geq 2$, there exists a
constant $c_k$ such that for all constants $c<c_k$,
$F(n,cn)$ is a.a.s.\ (asymptotically almost surely) satisfiable,
while for $c>c_k$ it is a.a.s.\ unsatisfiable.

The case of 2-SAT is well understood, with
Chv\'{a}tal and Reed~\cite{mick},
Geordt~\cite{GoTU}, and Fernandez de la Vega~\cite{vega2sat}
independently proving
that $c_2= 1$, and \bol, Borgs, Chayes, Kim, and Wilson
\cite{2satwin} determined the ``scaling window'' to be
$1+\Theta(n^{-1/3})$.
For $k> 2$, the conjecture remains open.
Friedgut proved that for any $k$ and $n$ there is sharp threshold $c_k(n)$,
leaving open whether $c_k(n)$ has a limit $c_k$.
With the threshold behavior not understood, considerable attention has
been devoted to proving density bounds below which a formula is a.a.s.\
satisfiable (``lower bounds'' on the putative threshold)
and bounds above which it is a.a.s.\ unsatisfiable (``upper bounds'').
For $k=3$, it is conjectured that $c_3 \approx 4.2$,
and the best upper bound is 4.596 \cite{Janson}.%
\footnote{A bound of $4.506$ due to
Dubois, Boufkhad, and Mandler~\cite{dub_announ}
has not appeared in journal-refereed form.}

Existing lower bounds for 3-SAT are all algorithmically based.
(By contrast, new lower bounds for $k$-SAT
are based on the second-moment method \cite{AM02,AP03}.)
The earliest such bound
of $1.63$, due to Broder, Frieze and Upfal~\cite{BFU}
was based on the ``pure-literal rule'':
successively setting to True literals whose complement does not appear,
``reducing'' the formula, and repeating.
The next, $3.003$, due to Frieze and Suen~\cite{FrSu},
used the ``shortest-clause rule'',
setting True a random literal from a random shortest clause.
Skipping over a bound of $3.145$ by Achlioptas~\cite{pi},
a bound of $3.14$ due to Achlioptas and Sorkin~\cite{AchSo00}
again selects a literal from a shortest clause,
but when the literal is from a 2-clause (the case of interest)
it sets the literal either
True or False depending (optimally) on the number of other occurrences
of the literal and its negation in 2- and 3-clauses;
\cite{AchSo00} extends this to a version which optimally
chooses to set one or two literals at a time, and sets them optimally,
for a bound of 3.26.
\cite{AchSo00} suggests that better bounds may require looking at literal-degree
information, in some way harking back to~\cite{BFU}.
This approach was taken up by Kaporis, Kirousis, and Lalas~\cite{KKL02},
whose algorithm sets a variable of largest degree
(a ``1-parameter heuristic'')
to give a bound of 3.42.
It was clear that the same approach could be exploited further.

\section{Result, significance, and open problems}
In this paper, we choose a variable according to its degree and that of
its complement (a 2-parameter heuristic), to get a bound of 3.52.
Kaporis and Lalas~\cite{KL},
using a similar but not identical heuristic,
independently obtained the same bound at around the same time.
The purpose of this short abstract is twofold.

First, since the bounds from all heuristics of this sort rely on
numerical calculations (notably, solutions to differential equations),
it is important to put on record that our calculations and those of
\cite{KL} independently justify a value of 3.52,
and that we reproduce the 3.42 bound of~\cite{KKL02}.

Second, our heuristic (and that of \cite{KL}) efficiently solves
denser random instances than any other theoretically justified algorithms;
since solving 3-SAT instances is of practical importance,
our algorithm may be of practical utility.
In that regard, a few remarks.
Our heuristic succeeds only with a probability that is asymptotically
bounded away from~0, but exploiting a standard one-step backtracking
trick brings the asymptotic probability to~1.
Also, the algorithms commonly used in practice are Davis-Putnam-type
backtracking procedures, quite different from the ``greedy'' approaches
taken in all the works described above.
However, it is easy to imagine
using the present heuristic as a selection rule for a Davis-Putnam
algorithm, preserving the heuristic's theoretically justified behavior on
random instances, while gaining the Davis-Putnam algorithm's guarantee
of a correct answer on arbitrary instances.

Not present in this short abstract is a rigorous justification of
our proof methodology (fairly easy and familiar), nor of the numerical
calculations, which to be done rigorously would require theoretically
derived Lipschitz bounds on a derivative, and interval-arithmetic
calculations employing those bounds, along with a few other
technicalities.

Future work could include consideration of a variable's number of
appearances and that of its complement separately in 2-clauses and
3-clauses (a 4-parameter heuristic),
which is analyzable in the same framework.
In at least the 2- and 4-parameter versions of literal-degree heuristics
(as opposed to the 1-parameter version),
it is not clear how best to select a next literal: does a $(2,3)$
literal (2 positive appearances, 3 negative) trump a $(4,5)$, or vice-versa?
In the 4-parameter version, it is also not clear how best to set a chosen
literal; this was the question answered
in~\cite{AchSo00}
for the non-degree-spectrum case.
An optimal solution to these questions would be a most interesting
theoretical contribution, and could also give significant
improvements in the bounds.

\section{Algorithm}
\label{sec1}
We call a variable with $i$ positive and $j$ negative appearances
an $(i,j)$-variable.
Our algorithm is defined as follows.

\vspace{-0.3cm}

{\small
\label{alpage}
\begin{tabbing}
1233\=1233\=4566\=7899\=78999\=7899\=789\=789\=789\=789\=\kill
{\bf Algorithm A}\\
{\em Input:} \>\> A 3-CNF formula. \\
{\bf begin}\\
1 \> {\bf while} there exists an unset variable\\ 
2 \>\> {\bf choose}  an $(i,j)$-degree variable  using a selection rule\\
3 \>\> {\bf set} $v$ True if $i< j$ and False otherwise\\
4 \>\> {\bf while} there exists a unit clause \\
5 \>\>\> {\bf set} a literal of an arbitrary unit clause True.\\
6 \> {\bf if} an empty clause is generated report {\em failure};
otherwise report {\em success}\\
{\bf end}\\
\end{tabbing}
} \vspace{-0.75cm}
The best selection rule we found was this.
If there is a ``pure'' variable (one with $i=0$ or $j=0$), select it.
Otherwise, choose a variable with maximum discrepancy $|i-j|$,
breaking ties in favor of maximal $i+j$.
(This identifies a unique unordered pair $\{i,j\}$,
and all variables with those degrees are indistinguishable
to the algorithm.)
This selection rule satisfies formulas up to density~3.52.
Other selection rules we tried were less good.
Working as above but breaking ties in favor of \emph{minimal} $i+j$
only worked up to density~3.50.
Selecting by maximum $i/j$ instead of $i-j$ only worked up to~3.44.
Selecting by maximum $\max\{i,j\}$ is equivalent to the approach of
\cite{KKL02}, and we reproduce their~3.42.

\section{Analysis}
In truth, the ``natural'' algorithm above is not the one analyzed.
Rather than making $\Theta(n)$ iterations, the analyzed algorithm
makes a constant number of iterations, in each of which
it sets $\Theta(n)$ variables with common degree $\{i,j\}$.

It is easily verified that during the algorithm, the formula remains
uniformly random conditioned on its degree sequence.
To make the calculations finite,
we truncate the degree sequence at
some value $h$ ($h=31$ in our calculations).
Then with $n$ the original number of variables,
for $i,j < h$
we let $n_{i,j}$ be $1/n$th the current number of variables of degree $(i,j)$;
$n_{h,j}$ (and $n_{i,h}$) the similar value for variables
with $\geq h$ positive (negative) appearances;
and $n_{h,h}$ that for variables with $\geq h$ positive and negative
appearances.
Setting a single $(i,j)$-variable produces straightforwardly computable
expected changes $\Delta$
(detailed in Appendix~A) to
the $h^2$-dimensional vector $S$ of values $n_{i,j}$,
and each element of $\Delta$ has order only $O(1/n)$,
so the differential equation method
(see for example \cite{worm})
can be used to prove that, as we set $\Theta(n)$ variables
with common degree $(i,j)$,
the vector $n_{i,j}$
almost surely almost exactly follows a trajectory described by
the solution of a differential equation corresponding to~$\Delta$.

So instead of selecting an $(i,j)$-variable as in Algorithm~A,
we use the same selection rule to select a pair
$(i,j)$, $i,j \leq h$,
and we set $n \cdot \min\{\delta, n_{i,j}\}$ $(i,j)$-variables at once.
Here $\delta$ is a value of our choosing, which could vary from round
to round, but which we fixed at $10^{-6}$.
Each such round (including the unit-clause steps it implies)
can be described by the differential equation method,
and our analysis simply consists of simulating the differential equations
for a constant number of rounds.
It is clear that after some number of rounds, all the values in $S$
can be made arbitrarily small,
and at that point we apply the main theorem of
Cooper, Frieze, and Sorkin \cite{CFS} to show that the remaining formula
is satisfiable a.a.s.
(The positive side of their result has a natural algorithmic interpretation,
so our procedure remains algorithmic to the end.)

\section{Differential Equations}
\label{appendix} In this section, we describe the
differential equations for the case in which we have a
2-dimensional table for keeping the expected number of variables with $k<
h$ positive appearances and $l< h$ negative appearances. Since
$dt$ (the ``time parameter'' described before) is very small,
w.l.o.g.\ we can assume all values of $S$ remain fixed during a
round. Using these values, we obtain the new value of $S$ after a
round. Suppose we set a variable from cell $n_{k,l}$ True (in other
words, set a $(k,l)$-degree variable True).
Then writing $\delta$ to denote the \emph{expected} increase
to a parameter, for such a ``free move'' we have:
\begin{align*}
\delta m_3 &= -k\frac{3 \cdot m_3}{\Sigma_m}- l\frac{3\cdot
m_3}{\Sigma_m} ,
\\
\delta m_2 &= -k\frac{2 \cdot m_2}{\Sigma_m}-
l\frac{2\cdot m_2}{\Sigma_m}+
 l\frac{3 \cdot m_3}{\Sigma_m} ,
\\
\delta m_1 &= l\frac{2 \cdot m_2}{\Sigma_m} , \text{ and}
\\
\delta n_{ij} &= -\psi(i,k)\psi(j,l)
 -k \cdot \frac{2\cdot 1\cdot m_2+ 3\cdot2\cdot m_3}{\Sigma_m} \cdot
 \\ & \qquad
 \cdot \frac{(i+j)\cdot n_{ij}-(i+1)\cdot n_{i+1,j}- (j+1)\cdot
n_{i,j+1}}{\Sigma_n},
\end{align*}
where $\Sigma_m=\Sigma_n=2\cdot m_2+3
\cdot m_3$ is the total density of appearances of all variables
and $m_1$ is the expected number of unit clauses generated by this free
move. Here $\psi(x,y)= 1$ if $x=y$ and zero otherwise.
Note also that since by definition $m_1=0$ at the start of a round,
at the end, $m_1 = \delta m_1$ as given above.

After a free move, we have a number of ``forced moves'' in which
the literals in all $m_1$ unit clauses must be set  True to satisfy our
formula.
 A literal in a
unit
clause is a variable from cell $({k+1,l})$
with probability $\frac{(k+1)n_{k+1, l}}{\Sigma_n}$,
or the negation
of a variable from cell $({k,l+ 1})$ with probability
$\frac{(l+1)n_{k,l+1}}{\Sigma_n}$.
In either case, in the rest of the formula (excepting the unit clauses)
that variable has degree $(k,l)$. Thus the expected number $\rho$ of
new unit clauses produced by
one such forced move
(the Malthus parameter in our Galton-Watson process)
is
\begin{align*}
\sum_{0\leq k', l'<
h}\frac{(k'+1)n_{k'+1, l'}}{\Sigma_n}l'\frac{ 2 \cdot
m_2}{\Sigma_m}+ \frac{(l'+1)n_{k', l'+1}}{\Sigma_n}k' \frac{ 2
\cdot m_2}{\Sigma_m} ,
\end{align*}
where $l'\frac{ 2 \cdot m_2}{\Sigma_m}$
(see parameter $m_1$ defined above) is the density of new unit
clauses after setting a $(k',l')$-degree variable True (which
happens with probability $\frac{(k'+1)n_{k'+1, l'}}{\Sigma_n}$) and  $k' \frac{ 2 \cdot m_2}{\Sigma_m}$ is the
density of new unit clauses after setting a $(k',l')$-degree
variable False (which happens with probability
$\frac{(l'+1)n_{k', l'+1}}{\Sigma_n}$). For such a forced move,
the expected parameter changes are:
\begin{align*}
\delta' m_3 &=
\sum_{0\leq k', l'< h}\frac{(k'+1)n_{k'+1, l'}}{\Sigma_n}\delta
M_3(k', l', {\text T}) + \frac{(l'+1)n_{k', l'+1}}{\Sigma_n}\delta
M_3(k', l', {\text F}) , 
\\
\delta' m_2 &= \sum_{0\leq k', l'<
h}\frac{(k'+1)n_{k'+1, l'}}{\Sigma_n}\delta M_2(k', l', {\text T})
+ \frac{(l'+1)n_{k', l'+1}}{\Sigma_n}\delta M_2(k', l', {\text
F}) ,
\end{align*}
where $\delta M_3(k', l', {\text T})$
has exactly the same formula as $\delta m_3$ defined above,
likewise for $\delta M_2(k', l', {\text T})$ and
$\delta m_2$,
and where by symmetry $\delta M_3(k', l', {\text F})=\delta
M_3(l', k', {\text T})$
and $\delta M_2(k', l', {\text F})=\delta M_2(l',k', {\text T})$.

Finally, for each $i$ and $j$,
\begin{align*}
\delta' n_{i,j} &=  \sum_{0\leq k', l'< h} \frac{(k'+1)n_{k'+1,
l'}}{\Sigma_n}(\delta
N_{i,j}(k', l', {\text T})
 \\ & \qquad
 -\psi(i,k'+1)\cdot \psi(j,l'))
 \\ & \qquad
+ \frac{(l'+1)n_{k',
l'+1}}{\Sigma_n}(\delta N_{i,j}(k', l', {\text F})- \psi(i,k')\cdot
\psi(j,l'+1)) ,
\\ \intertext{where}
\delta N_{i,j}(k, l, {\text T})
 &= -k \cdot \frac{2 \cdot 1 \cdot m_2+
   3\cdot 2\cdot m_3}{\Sigma_m} \; \cdot
   \\ & \qquad
   \cdot \left( \frac{(i+j)\cdot n_{i,j}}{\Sigma_n}-\frac{(i+1)
   \cdot n_{i+1,j}}{\Sigma_n}-\frac{(j+1)\cdot n_{i, j+1}}{\Sigma_n} \right)
  \text{ and}
\\
\delta N_{i,j}(k, l, {\text F})
 &= \delta N_{i,j}( l, k, {\text T}) .
\end{align*}
We note that the formula for $\delta' n_{i,j}$ can be
obtained by considering the flow which goes in or out for cell $n_{i,j}$.

\begin{sloppy}
Reasoning via the
Galton-Watson process, we know that 
the expected  number of forced moves is $\frac{m_1}{1-\rho}$.
Thus the new expected value of $S$ after setting a small fraction $dt$ of
variables from cell $n_{k,l}$  True is: $S+ dt(\delta m_2, \delta
m_3, \delta n)+ dt\frac{m_1}{1-\rho}(\delta' m_2, \delta'
m_3,\delta' n)$. If we set a variable
from cell $n_{k,l}$  False, the expected changes can be obtained by just
swapping the role of $k$ and $l$ in the above description.
\end{sloppy}

\section*{Acknowledgments}
The differential equations we used were derived by Frieze and Sorkin
in 2001, and we are grateful to Alan Frieze for his ongoing help.
We also appreciate open exchanges with
Alexis Kaporis, Lefteris Kirousis, and Efthimios Lalas.

\newcommand{\etalchar}[1]{$^{#1}$}
\providecommand{\bysame}{\leavevmode\hbox to3em{\hrulefill}\thinspace}

\end{document}